\documentclass[12pt]{article}
\usepackage{graphicx,bm}
\renewcommand\[{\begin{equation}}
\renewcommand\]{\end{equation}}

\begin{document}

\begin{center}
{\large \bf Rusanov-type Schemes for Hyperbolic Equations: \\
Wave-Speed Estimates, Monotonicity and Stability}
\end{center}

\begin{center} \vspace*{15pt} 
{\bf Eleuterio F.  Toro~$^a$ and Svetlana A. Tokareva~$^{b}$}
\end{center}

\begin{center}
$^{a}$ Laboratory of Applied Mathematics, DICAM, University of Trento, Italy \\
Email: eleuterio.toro@unitn.it
\end{center}


\begin{center}
$^{b}$ Los Alamos National Laboratory, Theoretical Division, MS B284, 
Los Alamos, NM, 87545, USA.  \\
Email: tokareva@lanl.gov
\end{center}

{\it In memory of  our long-standing collaborator and friend Arturo Hidalgo,  Universidad Polit\'ecnica de Madrid, Spain.}\\

\begin{abstract}

HLL-type schemes constitute a large hierarchy of numerical methods,  in the finite volume and discontinuous Galerkin finite element frameworks, for solving hyperbolic equations. The hierarchy of fluxes includes Rusanov schemes, HLL schemes,  HLLC schemes,  and other variations.  All of these schemes rely on wave speed estimates.  Recent work 
\cite{Toro:2020c} has shown that most  wave speed estimates in current use underestimate the true wave speeds.  In the present paper we carry out a theoretical study of the consequences arising from errors in the wave speed estimates,  on the monotonicity and stability properties of the derived schemes.   For the simplest case of the hierarchy,  that is Rusanov-type schemes,  we carry out  a detailed analysis in terms of the  linear advection equation in one and two space dimensions.   It is found that errors from underestimates of the wave speed could cause  loss of monotonicity,  severe reduction of the stability limit,  and even loss of stability.   Errors from overestimates,  though preserving monotonicity,  will cause a reduction of the stability limit.  We find that overestimation is  preferable to underestimation,  for two reasons.  First,  schemes from overestimation are monotone,  and second, their stability regions are larger than those from underestimation,  for equivalent displacements from the exact speed.  The findings of this paper may prove useful in raising awareness of the potential pitfalls of a seemingly simple practical computational task, that of providing wave speed estimates. Our reported findings may also motivate subsequent studies for complex non-linear hyperbolic systems,  requiring estimates for two or more waves,  such as in HLL and HLLC schemes.

\end{abstract}


\section{Introduction}

The Rusanov flux  \cite{Rusanov:1961a},  often called  the {\it local  Lax-Friedrichs} flux \cite{Leveque:2002a},  is  frequently used in finite volume and discontinuous Galerkin finite element methods for solving hyperbolic equations 
\cite{Shu:1987b},
\cite{Titarev:2004c}, 
\cite{Bouchut:2004a},
\cite{Cheng:2007a}, 
\cite{Toro:2009a},
\cite{CastroM:2012a},
\cite{Dumbser:2016d},
\cite{Fambri:2018a},
\cite{Rannabauer:2018a},
\cite{Busto:2020a},
\cite{Kuzmin:2021a},
\cite{Gaburro:2021a},
\cite{Mohamed:2021a},
\cite{Busto:2021a},
\cite{Gomez-Bueno:2023a},
\cite{Abgrall:2023a}.
The Rusanov method  is  the simplest upwind scheme of the Godunov-type  \cite{Godunov:1959a},  requiring a single wave-speed estimate $\hat s$ to fully determine the flux.  That is,  the Rusanov scheme adopts a one-wave model,  as it accounts for just one wave,  the fastest from the wave configurations triggered at cell interfaces by the associated Riemann problems \cite{Toro:2009a}.  The more sophisticated HLL scheme \cite{Harten:1983b} introduces a two-wave model,  thus requiring estimates for the minimal $S_{L}$ and maximal $S_{R}$ wave speeds.  The HLLC scheme \cite{Toro:1994c} adopts a three-wave model,  thus requiring estimates for three  wave speeds.  The Godunov method in conjunction with the exact Riemann problem accounts for all $m$ waves present in the solution of the Riemann problem for a hyperbolic system of $m$ equations.  Schemes adopting a  z-wave model, with $z<m$ are said to be {\bf incomplete},  while those  with $z=m$ are said to be {\bf complete}.  The Rusanov scheme is incomplete for any system of equations ($m>1$).  Note that all incomplete fluxes sacrifice the resolution of intermediate characteristic fields,  the effect being most obvious for the case of slowly moving,  or stationary,  intermediate waves.

The focus of this paper  concerns the one-wave model Rusanov scheme.   In general,  the sought exact wave speed $\hat s$ is unknown and it may be very costly,  or even impossible,  to compute.  Therefore,  the  Rusanov scheme relies on estimates $\hat s$,  which from their approximate nature,  will invariably carry an error.   The topic of wave speed estimates for hyperbolic equations, in a computational setting,  is more substantial than given credit for in the literature.  Reliable estimates, for example,  are needed for enforcing the Courant stability condition and choose a stable time step in explicit methods; they may also be used in constructing numerical fluxes of the HLL-type \cite{Harten:1983b},  as already pointed out.   Front-tracking and shock-fitting methods also need accurate wave speed estimates. 
And yet, most existing methods for estimating wave speeds are inaccurate and fail to bound exact wave speeds 
\cite{Guermond:2016a}, \cite{Toro:2020c}.  As necessary background for this paper,  we first  address the problem of estimating wave speeds and review most methods available for doing so; we also illustrates the challenges in obtaining bounds and  discuss potential applications of wave speed estimates in a computational setting.

The main subject of this paper regards a detailed analysis of the consequences from erroneous estimation on the monotonicity and stability properties of the derived numerical schemes.  The study is performed in terms of the model linear advection equation in one and two space dimensions on Cartesian meshes.  
For the one-dimensional case the exact speed is denoted by $\lambda$ and it is assumed that the uncertainty arising from estimating $\hat s$  can be represented by a parameter  $\beta$,  such that $\hat s= \beta \lambda$,  with $\beta \in [1-\epsilon_{B}, 1+\epsilon_{T}]$,  $0\le \epsilon_{B} \le 1$ and $0 \le \epsilon_{T} < \infty$.   In one space dimension the exact wave speed is recovered when $\beta=1$,  thus identically reproducing the Godunov upwind method.
We show that the estimate $\hat s(\beta)$ reproduces a broad family of well-known schemes,  such as the Godunov upwind method \cite{Godunov:1959a};  the Lax-Wendroff method  \cite{Lax:1960a};  the Godunov centred method \cite{Godunov:1962a};  the FORCE method \cite{Toro:2000a};  the  FORCE-$\alpha$ method \cite{Toro:2020b} and the Lax-Friedrichs method \cite{Lax:1960a},  \cite{Crandall:1980a}.  Our analysis establishes that the choice of the wave speed estimate 
$\hat s$ has a profound effect on the monotonicity and linear stability properties of the resulting numerical methods in one and two space dimensions.  Some of the results will extrapolate to a wider class of schemes,  including the 
HLL and HLLC  numerical methods for approximating systems of non-linear hyperbolic equations.

The rest of the paper is structured as follows.  Section \ref{sec:Preliminaries} sets the strictly necessary background and motivation for this work;  Section \ref{sec:Rusanov} formulates Rusanov-type schemes for the one-dimensional linear advection equation; Section \ref{sec:Perturbation} discusses some of the consequences of perturbing the exact speed on monotonicity and stability of the family of Rusanov-type schemes; Section \ref{sec:Rusanov2D} extends the analysis to Rusanov-type schemes as applied  to the linear advection equation in two space dimensions; a summary and concluding remarks are found in Section \ref{sec:Conclusions}.

\section{Preliminaries: Wave Speeds and Rusanov Fluxes}
\label{sec:Preliminaries}

In this section we provide the strictly necessary background for discussing the main theme of this paper developed in Sec.
\ref{sec:Rusanov},  which is Rusanov-type schemes that emerge from estimates for the single wave speed $\hat s$ needed to define the numerical flux.   Here we address the problem of estimating wave speeds,  illustrates the challenges in obtaining bounds and  discuss potential applications of wave speed estimates in a computational setting.

Let us consider the initial-boundary value problem (IBVP)  for a system of $m$ hyperbolic conservation laws
\begin{equation}                                                  \label{eq:PREL1}
      \left.\begin{array}{ll}
       \mbox{PDEs: } &  \partial_{t}{\bf Q} + \partial_{x}{\bf F}({\bf Q}) = {\bf 0} \;,\hspace{2mm} x \in [a,b] \;, \hspace{2mm} t>0 \;, \\
       \mbox{ICs: }  &  {\bf Q}(x,0) = {\bf Q}^{(0)}(x)\;, \hspace{2mm} x \in [a,b] \;, \\
       \mbox{BCs: }  &  {\bf Q}(a,t) = {\bf B}_{L}(t)\;, \hspace{2mm}{\bf Q}(b,t) = {\bf B}_{R}(t)\;, \hspace{2mm} t \ge 0 \;.
      \end{array}\right\}
\end{equation}
Here ${\bf Q}(x,t)$ is the vector of unknown conserved variables;  ${\bf F}({\bf Q})$ is the physical flux vector;  ${\bf Q}^{(0)}(x)$ is the initial condition; ${\bf B}_{L}(t)$ and ${\bf B}_{R}(t)$ define boundary conditions.  Consider  also the  conservative method 
\begin{equation}                                                \label{eq:PREL2}
      {\bf Q}_{i}^{n+1} =  
      {\bf Q}_{i}^{n} -\frac{\Delta t}{\Delta x} \left( {\bf F}_{i+\frac{1}{2}} -{\bf F}_{i-\frac{1}{2}} \right) \;
\end{equation}
to solve (\ref{eq:PREL1}) on a regular mesh of size $\Delta x$, with time step $\Delta t$ computed on stability grounds. In general, the numerical flux is
\begin{equation}	                                             \label{eq:PREL3}
    {\bf F}_{i+\frac{1}{2}} = \frac{1}{\Delta t} \int^{\Delta t}_{0} {\bf F} ({\bf Q}(x_{i+\frac{1}{2}},t)) dt  \;.
\end{equation}
To determine ${\bf Q}(x_{i+\frac{1}{2}},t)$,  and hence the numerical flux from (\ref{eq:PREL3}),  Godunov suggested \cite{Godunov:1959a} to solve the Riemann problem
\begin{equation}                                               \label{eq:PREL4}
	    \left. \begin{array}{ll}
		  \mbox{PDEs:} & \partial_{t}{\bf Q} + \partial_{x} {\bf F}({\bf Q})={\bf 0} \;,  \\
		  \mbox{ICs:}  & {\bf Q} (x,0) = 
			\left\{ \begin{array}{lll}
			 {\bf Q}_{L} = {\bf Q}_{i}^{n}   &  \mbox{ if } & x < x_{i+\frac{1}{2}}   \;,  \\
		    {\bf Q}_{R} = {\bf Q}_{i+1}^{n} &  \mbox{ if } & x >x_{i+\frac{1}{2}}   \;. 
			\end{array} \right.
	    \end{array} \right\}
\end{equation}
The similarity solution ${\bf Q}_{i+\frac{1}{2}}(x/t)$ of (\ref{eq:PREL4}) will exhibit $m$ wave families associated with $m$ real eigenvalues.
${\bf Q}_{i+\frac{1}{2}}(x/t)$ is found in local coordinates and evaluated at the interface $x/t=0$,  thus giving the Godunov state ${\bf Q}_{i+\frac{1}{2}}(0)$.  Noting that ${\bf Q}_{i+\frac{1}{2}}(0)$ is constant along the $t$-axis,  the Godunov flux emerges from  (\ref{eq:PREL3}) as
\begin{equation}                                              \label{eq:PREL5}
      {\bf F}_{i+\frac{1}{2}}={\bf F}({\bf Q}_{i+\frac{1}{2}}(0)) \;.
\end{equation} 
Detailed background on Godunov methods and solvers for the Riemann problem are found in \cite{Toro:2001a}, \cite{Toro:2009a}, \cite{Toro:2024a}. 

\subsection{Wave speeds for computing the time step}

Explicit schemes of the form (\ref{eq:PREL2}) require the determination of a suitable time step $\Delta t$.  This is accomplished by enforcing the Courant-Friedrichs-Lewy (CFL) condition, usually based on a linear stability analysis of the scheme for the model linear advection equation 
\begin{equation}                                              \label{eq:CFL-1}
\partial_{t}q +  \lambda \partial_{x}q=0 \;, 
\end{equation}
where $\lambda$ is the (constant) wave speed.  Assuming the stability limit to be $c_{lim}$, one requires
\begin{equation}                                              \label{eq:CFL-2}
      c= \frac{\Delta t \lambda}{\Delta x} \le c_{lim} \hspace{3mm}  \rightarrow  \hspace{3mm} \Delta t \le c_{lim} \frac{\Delta x}{\lambda} \;, 
\end{equation} 
from which we write
\begin{equation}                                              \label{eq:CFL-3}
        \Delta t = C_{cfl} \times c_{lim} \frac{\Delta x}{\lambda} \;.
\end{equation} 
Here $c$ is the CFL number and $C_{cfl}$ is the CFL (safety) coefficient, with $0<C_{cfl} \le 1$.   Most well known methods in one space dimension have $c_{lim}=1$.  One exception is Godunov's centred method  \cite{Godunov:1962a}, which has $c_{lim}=\frac{1}{2}\sqrt{2}$.

 For a nonlinear system,  such as the 1D Euler equations,  one generalises (\ref{eq:CFL-3}) as
\begin{equation}                                              \label{eq:CFL-4}
        \Delta t = C_{cfl} \times c_{lim} \frac{\Delta x}{S_{max}^{n}} \;.
\end{equation} 
Computing $ \Delta t$ from (\ref{eq:CFL-4}) requires a suitable  estimate for the wave speed $S_{max}^{n}$ at time level $n$.  In practice, almost invariably, one uses the eigenvalues of the system,  namely $\lambda_{1}=u-c$, $\lambda_{2}=u$, $\lambda_{3}=u+c$, where $u$ is the particle velocity and $c$ is the speed of sound. Typically, these are evaluated at the data states of the Riemann problem   (\ref{eq:PREL4}),  implying the choice
\begin{equation}                                              \label{eq:CFL-5}
       S_{max}^{n}= \max_{i} \left\{  |u_{i}^{n}| + c_{i}^{n} \right\} \;.
\end{equation} 
It is obvious that this choice fails to bound all possible wave speeds arising from (\ref{eq:PREL4}).  Consider for example the case in which at the initial time the particle velocity is zero everywhere.  In this case one hopes to bound shock speeds purely via the sound speed evaluated on the data states.  Choice (\ref{eq:CFL-5}) is bound to  result in an underestimate of the wave speed, leading to  an overestimate of the time step $ \Delta t$ in (\ref{eq:CFL-4}),  possibly rendering the scheme unstable at the beginning of the calculations.  Of course one can use the safety coefficient $C_{cfl}$ in (\ref{eq:CFL-4}) to attempt to remedy the situation, in a largely, and sometimes frustrating,  empirical fashion.  If bounds for the maximal true wave speed were available, then one could confidently take $C_{cfl}=1$, or just below unity to cater for round off errors.  Note that acceleration of waves resulting from nonlinear wave interaction  within a cell in time $\Delta t$  is possible,  which may require a net reduction of $C_{cfl}$ below unity.
 
\subsection{Wave speeds for computing fluxes}
\label{sec:WSE}

Bounds for the true wave speeds are also useful in the family of HLL-type numerical fluxes. Harten, Lax and van Leer \cite{Harten:1983b} proposed an approximation the solution of the Riemann problem (\ref{eq:PREL4}) to define a numerical flux for
(\ref{eq:PREL2}), called HLL flux.  Unlike the Godunov flux (\ref{eq:PREL5}) that first computes a state ${\bf Q}_{i+\frac{1}{2}}(0)$,  HLL gives a direct approximation to the flux  ${\bf F}_{i+\frac{1}{2}}$ in (\ref{eq:PREL3}),  as follows:
\begin{equation}	                                      \label{eq:PREL6}
  {\bf F}^{HLL}_{i+\frac{1}{2}}=
  \left\{\begin{array}{cc}
  {\bf F}({\bf Q}_{i}^{n})                                                                                                &  0 \leq S_{L} \;,  \\
 \displaystyle{\frac{S_{R}{\bf F}({\bf Q}_{i}^{n})-S_{L}{\bf F}( {\bf Q}_{i+1}^{n})+S_{L}S_{R}({\bf Q}_{i+1}^{n}-
  {\bf Q}_{i}^{n})}{S_{R}- S_{L}}}\;, &  S_{L} < 0 <S_{R} \;, \\
  {\bf F}({\bf Q}_{i+1}^{n})                                                                                             & 0 \geq S_{R} \;.
  \end{array}\right.
\end{equation}
The HLL flux (\ref{eq:PREL6}) is fully determined once estimates $S_{L}$ and $S_{R}$ for the minimal and maximal wave speeds present in the solution of (\ref{eq:PREL4}) have been provided.  More details on HLL are found in  \cite{Harten:1983b}, \cite{Toro:2001a}, \cite{Toro:2009a}, \cite{Toro:2024a}.  For an up to date discussion on wave speed estimates see 
\cite{Guermond:2016a}, \cite{Toro:2020c} and references therein.  In what follows we list most of the well-known wave speed estimates available in the literature. We do so for the one-dimensional ideal Euler equations.

\begin{enumerate}
\item Davis \cite{Davis:1988a} proposed two ways to estimate the speeds $S_{L}$ and $S_{R}$, namely
\begin{equation}           \label{Davis:a}
       S^{Dav_a}_{L}= u_{L} - c_{L} \;, \qquad S^{Dav_a}_{R}= u_{R} + c_{R} \;
\end{equation}
and 
\begin{equation}           \label{Davis:b}
         S^{Dav_b}_{L}= \min\{u_{L} - c_{L} , u_{R} - c_{R}\} \;, \hspace{1mm} S^{Dav_b}_{R}= \max\{u_{L} + c_{L} , u_{R} + c_{R}\}\;.
\end{equation}

\item Einfeldt  \cite{Einfeldt:1988a} proposed 
\begin{equation}
       S^{Einf}_L= \widetilde{u} - \widetilde{d} \;,  \qquad S^{Einf}_R= \widetilde{u} + \widetilde{d}   \;,
\end{equation}
where 
\begin{equation} \label{EQ:utilde}
      \widetilde{u}= \frac{\sqrt{\rho_L}u_L + \sqrt{\rho_R}u_R}{\sqrt{\rho_L} + \sqrt{\rho_R}} \qquad 
\end{equation}
and
\begin{equation} \label{EQ:eta2}
      {\widetilde{d^2}} =\frac{\sqrt{\rho_L}c^2_L + \sqrt{\rho_R}c^2_R}{\sqrt{\rho_L} + \sqrt{\rho_R}} +\frac{1}{2} \frac{\sqrt{\rho_L} \sqrt{\rho_R}}  
      {\left(\sqrt{\rho_L} + \sqrt{\rho_R}\right)^2} \left(u_R-u_L\right)^2 \;.
\end{equation}
Einfeldt's estimates  incorporate information from the Riemann problem solution, via the Roe averages (\ref{EQ:utilde}) from the approximate Roe Riemann solver for the ideal Euler equations \cite{Roe:1981a}.

\item Batten et al.  \cite{Batten:1997a} also suggested to use the Roe averages to obtain wave speed estimates for the Euler equations, namely
\begin{equation}                         \label{Battena}
     S^{Ba}_L= \min\{u_{L} - c_{L} , \widetilde{u} - \widetilde{c} \} \;, \qquad S^{Ba}_R= \max\{u_{R} + c_{R} , \widetilde{u} + \widetilde{c} \}  \;,
\end{equation}
where
\begin{equation}                         \label{Battenb}
     \left.\begin{array}{c}
     \displaystyle
     \widetilde{u}= \frac{\sqrt{\rho_L}u_L + \sqrt{\rho_R}u_R}{\sqrt{\rho_L} + \sqrt{\rho_R}} \;, 
     \widetilde{H}= \frac{\sqrt{\rho_L}H_L + \sqrt{\rho_R}H_R}{\sqrt{\rho_L} + \sqrt{\rho_R}} \;, \\
     \\
     \widetilde{c}= \left[(\gamma -1)\left(\widetilde{H}-\frac{1}{2}\widetilde{u}^2\right)\right]^{1/2} \;.
     \end{array}     \right\}
\end{equation}
$H_{K}$ is the specific enthalpy on the left and right states, namely
\begin{equation}                         \label{Battenc}
        H_{K} = \frac{E_{K}+p_{K}}{\rho_{K} } \quad(K=L,R) \;.
\end{equation}

\item Toro et al.   \cite{Toro:1992e}, \cite{Toro:1994c} proposed the estimates
\begin{equation}                          \label{Toroa}
    S_{L}^{To} = u_{L} - c_{L}q_{L}    \;, \hspace{3mm}  S_{R}^{To} = u_{R} +c_{R}q_{R} \;.
\end{equation}
The functions $q_{L}$ and  $q_{R}$, yet to be defined,  contain information from the Riemann problem solution and distinguish between shocks and rarefaction waves.   Successful choices are
\begin{equation}                          \label{Torob}
\begin{array}{c}
q_{L}=\left\{
       \begin{array}{cccc} 
       1  & \mbox{if} &  p_{*rr} \le p_{L}  & \mbox{left rarefaction} \;, \\
        \sqrt{1+\frac{\gamma+1}{2 \gamma}(\frac{p_{*rr}}{p_{L}}-1) }  & \mbox{if} &  p_{*rr} > p_{L}  & \mbox{left shock} 
       \end{array}\right.
\end{array}
\end{equation}
and
\begin{equation}                          \label{Toroc}
\begin{array}{c}
q_{R}=\left\{
       \begin{array}{cccc} 
       1  & \mbox{if} &  p_{*rr} \le p_{R}  & \mbox{right rarefaction} \;, \\
        \sqrt{1+\frac{\gamma+1}{2 \gamma}(\frac{p_{*rr}}{p_{R}}-1) }  & \mbox{if} &  p_{*rr} > p_{R}  & \mbox{right shock} \;.
       \end{array}\right.
\end{array}
\end{equation}
Here $p_{*rr}$ denotes an estimate for the pressure $p_{*}$ in the {\it Star Region} of the Riemann problem solution \cite{Toro:2009a}, such that 
\begin{equation}                          \label{Toroe}
       p_{*rr} \ge p_{*} \;.
\end{equation}
The following approximation fulfils condition (\ref{Toroe})
\begin{equation}                          \label{eeAllRarefaction}
      p_{ *rr}=\left[\frac{a_{ L}+a_{ R}-\frac{1}{2}(\gamma-1)(u_{ R}-u_{ L})}
      {a_{ L}/p_{ L}^{\frac{\gamma-1}{2\gamma}}
      +a_{ R}/p_{ R}^{\frac{\gamma-1}{2\gamma}}}
      \right]^{\frac{2\gamma}{\gamma-1}} \;.
\end{equation}
Choice (\ref{eeAllRarefaction}) used in (\ref{Toroa})-(\ref{Toroc}) provides theoretical bounds for the minimal and maximal intercell wave speeds \cite{Toro:2020c}.   
\end{enumerate}
Of the four choices listed, (\ref{Toroa}) are the only wave speed estimates that provide bounds for the true wave speeds for the Euler equations for ideal gases.  Such estimates can also be extended to the Euler equations with the co-volume equation of state \cite{Toro:1989d}.
For recent works on wave speed estimates see  Guermond and Popov \cite{Guermond:2016a}, \cite{Guermond:2024a} and Toro et al. \cite{Toro:2020c}.

\subsection{Numerical tests for wave speed estimates}
\label{sec:WSEexamples}

Here we reproduce some numerical experiments reported in Toro et al.  \cite{Toro:2020c}, where  seven Riemann problems for the ideal Euler equations with  $\gamma = 1.4$ were considered. Table \ref{tab:RPdata} gives the initial conditions  in columns 2 to 7; columns 8 and 9 give the exact solution for velocity $u_{*}$ and pressure  $p_{*}$ in the {\it Start Region},  while the last column describes the emerging wave pattern for each case.  
\begin{table}[h!]    
\begin{center}
\begin{tabular}{|c|c|c|c|c|c|c|c|c|c|} \hline
Test & $\rho_{L}$ & $u_{L}$  & $p_{L}$ & $\rho_{R}$ & $u_{R}$  & $p_{R}$ & $u_*$& $p_*$ & Waves  \\ \hline
1    &  1.0       & 0.0      & 1.0     &  1.0       &  0.0     & 0.1 &  0.5248    & 0.5219 & R-S    \\ \hline
2    &  1.0       & 0.0      & 1.0     &  0.125     &  0.0     & 0.1 &  0.9274 & 0.3031   & R-S \\ \hline
3    &  1.0       & 0.0      & 1.0     &  0.001     &  0.0     & 0.8 &  0.1794 & 0.8060    & R-S \\ \hline
4    &  1.0       & 0.0      & 0.01    &  1.0       &  0.0     & 1000.0 &  -19.5975 & 460.8938   & S-R \\ \hline
5    &  6.0       & 8.0      & 460.0   & 6.0        & -6.0     & 46.0&  3.8194  & 790.2928   & S-S \\ \hline 
6    &  600.0     & 80.0     & 4600.0  & 6.0        & -6.0     & 46.0&  44992.5781  & 790.2928   & S-S \\ \hline 
7    &  1.0       & -2.0     & 0.4     &  1.0       &  2.0     & 0.4 & 0.0000 &  0.0019    & R-R \\ \hline
\end{tabular}    
\caption{Initial conditions for seven  Riemann problems for the ideal Euler equations with $\gamma = 1.4$ and corresponding  exact solutions for  velocity $u_{*}$ and pressure $p_{*}$ in the {\it Star Region}.   Last column shows the emerging wave patterns, with S:shock and R:rarefaction.}
\label{tab:RPdata} 
\end{center}
\end{table}

Table \ref{Tab:results}   shows results for thr maximal wave speeds $S_R$ corresponding to the seven Rieman problems of Table \ref{tab:RPdata}. The exact solution $S_R^{Ex}$ is displayed in column 2, while current estimates discussed above are shown 
in columns 3 to 7.  It is seen that, with the exception of $S^{To}_{R}$ of Toro et al.   \cite{Toro:1992e}, \cite{Toro:1994c}, all remaining estimates fail to bound the maximal wave speed for all cases.  For recent advances on the subject of wave speed estimates, see  Guermond and Popov \cite{Guermond:2016a}, \cite{Guermond:2024a} and Toro et al. \cite{Toro:2020c}. 
\begin{table}[h!]
\begin{center}
\begin{tabular}{|c|c|c|c|c|c|c|} \hline
Test & $S^{Ex}_{R}$  & $S^{Dav_a}_{R}$ & $S^{Dav_b}_{R}$   & $S^{To}_{R}$  & $S^{Batten}_R$ &  $S^{Einf}_{R}$  \\ \hline
1  & 0.8039  & {\bf 0.3742} &          1.1832 &         0.8134    &      0.8775 &         0.8775  \\ \hline
2  & 1.7522 & {\bf 1.0583} &  {\bf 1.1832} &         1.7621 &  {\bf 1.1519} &     {\bf 1.1519}  \\ \hline
3  & 33.5742 & {\bf 33.4664} &    {\bf 33.4664} &        33.5742  &     {\bf 33.4664} &  {\bf 5.9740}  \\ \hline
4  & 37.4166 & 37.4166 &        37.4166 &        37.4166  &        37.4166 &  {\bf 26.4576} \\ \hline
5  & 6.6330 & {\bf -2.7238} &        18.3602 &         7.5400 &         9.2966 &        10.1397  \\ \hline
6  & 88.8686 & {\bf -2.7238} &  {\bf 83.2762} &       716.2437  &    {\bf 83.7136} &        89.9681  \\ \hline
7  & 2.7483 & 2.7483 &   2.7483 &   2.7483  &   2.7483 &  {\bf 1.6000}   \\ \hline
\end{tabular}
\caption{Maximal wave speed $S_R$  corresponding to the seven Rieman problems of Table \ref{tab:RPdata},  with exact solution  $S_R^{Ex}$ displayed in column 2.   Estimates are shown in columns 3 to 7.   Values in boldface fail to  bound the exact solution.}
\label{Tab:results}            
\end{center}
\end{table}

We have shown that most existing wave speed estimates  used in practice fail to bound the true wave speeds.   Moreover, the reported advances on the subject \cite{Guermond:2016a}, \cite{Guermond:2024a}, \cite{Toro:2020c} apply to the Euler equations for ideal gases.   Bound estimates  for compressible gas dynamics with general equations of state are yet to be found,  as are for other complicated hyperbolic systems.  For shallow water equations and for the blood flow equations the work is essentially complete \cite{Toro:2020c}.

\subsection{The Rusanov flux}
  
The {\bf Rusanov flux} \cite{Rusanov:1961a} is the simplest upwind flux that requires a wave speed estimate.  This flux can be derived as a special case of the HLL flux (\ref{eq:PREL6}),  for which wave speeds $S_{L}$ and $S_{R}$ are assumed to be available.  By first defining
\begin{equation}	                                      \label{eq:PREL7}
       \hat s =max\{|S_{L}|,|S_{R}|\} \;
\end{equation}
and then redefining $S_{L}$ and $S_{R}$ in (\ref{eq:PREL6}) as
\begin{equation}	                                      \label{eq:PREL7b}
       S_{R}=\hat s \;, \hspace{3mm} S_{L}= - \hat s \;, 
\end{equation}
the HLL flux (\ref{eq:PREL6}) reproduces the Rusanov flux  as follows
\begin{equation}	                                      \label{eq:PREL8}
   {\bf F}_{i+\frac{1}{2}} = \frac{1}{2} \left[ {\bf F}({\bf Q}_{i}^{n})  +  {\bf F}({\bf Q}_{i+1}^{n}) \right]
   -\frac{1}{2} \hat s ({\bf Q}_{i+1}^{n} -  {\bf Q}_{i}^{n} ) \;.
\end{equation}
To fully determine  (\ref{eq:PREL8}),  an estimate for the wave speed $\hat s$ must be prescribed, following (\ref{eq:PREL7}), for example,  applied to each intercell boundary,  at each time step.

Strictly speaking  \cite{Harten:1983b}, \cite{Toro:2020c},  the estimate $\hat s$ to be provided should bound  the true maximal wave speed in absolute value emerging from the exact solution of the Riemann problem (\ref{eq:PREL4}).  However,  as discussed in Secs.  \ref{sec:WSE} and \ref{sec:WSEexamples},  this task may prove difficult, costly or impossible to accomplish for sufficiently complex hyperbolic systems.  

In the next section we analyse the potential consequences resulting from inaccurate wave speed estimates.  We do so in terms of the simplest case,  namely the linear advection equation in one and two space dimensions. Obviously, for linear problems with constant coefficients the true wave speeds are known, which is advantageous for a precise,  quantitative assessment of the consequences of erroneous wave speed estimates.  We account for the uncertainties arising in non-linear systems by assuming wave speed estimates in a range that includes the true speeds in linear problems.

\section{Rusanov-Type Schemes}
\label{sec:Rusanov}

In this section we adopt a linear model equation in which the exact speed is known, thus permitting a quantitative,  theoretical assessment of the consequences derived from inaccurate wave speed estimates $\hat s$.

\subsection{Linear model}

\begin{figure}[h]
      \centerline{
             \includegraphics[scale=0.9, angle=0]{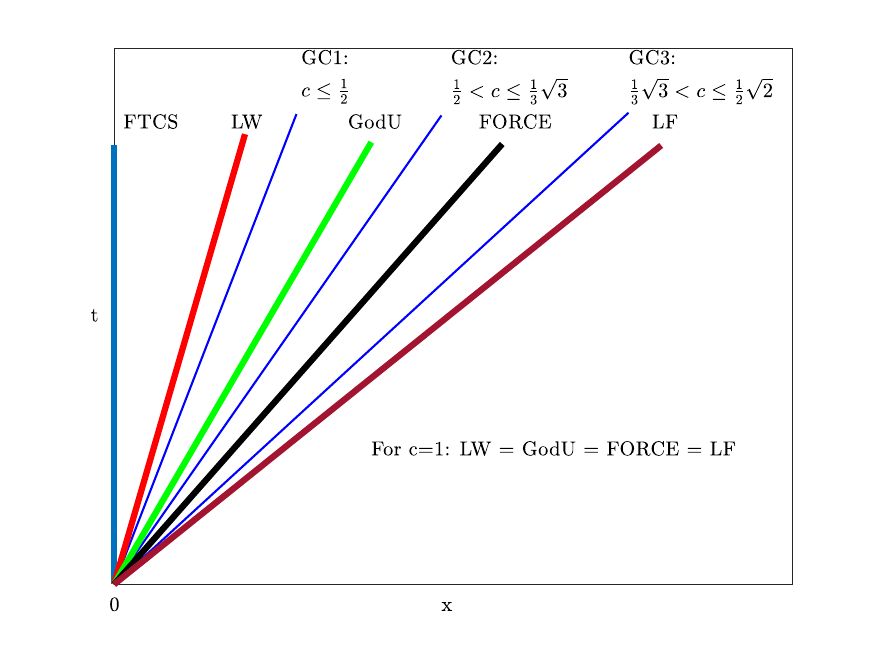}
      }
      \caption{Representation of Rusanov-type schemes in the $x$-$t$ plane.  Here FTCS is the Forward in Time Centred in Space scheme; LW is Lax-Wendroff; GodU is Godunov upwind; FORCE is the FORCE scheme and LF is Lax-Friedrichs. The three cases of the Godunov centred scheme are represented by GC1, GC2 and GC3. Each scheme corresponds to a particular choice of the characteristic line $\frac{x}{t}=\hat s$ emerging from the  origin $0$.  All linearly stable schemes lie in the wedge $LW0LF$. Large values of $\hat s$ are associated with monotone,  more diffusive schemes, while low values of $\hat s$ are associated with non-monotone schemes.}
      \label{fig:RusanovSchemes}
\end{figure}

In the present study we first consider the one-dimensional linear advection equation with constant and positive characteristic speed $\lambda$
\begin{equation}                                          \label{LAE-1}
    \partial_{t}q +  \partial_{x}f(q)= 0 \;, \hspace{3mm} f(q) = \lambda q \;, \hspace{3mm}  \lambda >0\;, 
\end{equation}
along with the conservative numerical scheme
\begin{equation}                                          \label{LAE-2}
    q_{i}^{n+1} = q_{i}^{n} - \frac{\Delta t}{\Delta x} [f_{i+\frac{1}{2}} - f_{i-\frac{1}{2}}] \;.
\end{equation}
The Rusanov \cite{Rusanov:1961a} numerical flux (\ref{eq:PREL8}) now reads
\begin{equation}                                          \label{LAE-3}
      f^{Rus}_{i+\frac{1}{2}} = \frac{1}{2}[f(q^{n}_{i}) +  f(q^{n}_{i+1})]
                                - \frac{1}{2} \hat s ( q^{n}_{i+1} - q^{n}_{i} ) \;, 
\end{equation}
where $\hat s$ is an estimate for the wave speed $\lambda$ of the differential equation (\ref{LAE-1}). 
Obviously, for equation (\ref{LAE-1}) the exact speed is known a priori, but for general systems this is not always the case,  and for which wave-speed estimates need to be proposed.  See \cite{Guermond:2016a}, \cite{Toro:2020c} for a discussion on wave speed estimates for systems.  In order to account for uncertainties in estimating the exact speed $\lambda$ by means of $\hat s$,  we shall assume
\begin{equation}                                    \label{LAE-4}
      \hat s (\beta)= \beta  \lambda  \;,
\end{equation}
with
\begin{equation}                                    \label{LAE-5}
      \beta \in [1-\epsilon_{B}, 1+\epsilon_{T}]\;,   \hspace{3mm} 0\le \epsilon_{B} \le 1 \;,  
       \hspace{3mm} 1 \le \epsilon_{T} < \infty \;.
\end{equation}
When $\beta=1$,  the exact wave speed is recovered,  thus identically reproducing the Godunov upwind method for (\ref{LAE-1}).  Other values of $\beta$ in (\ref{LAE-4}) reproduce different well-known schemes, as shown by some examples  in Table  \ref{Table:betaschemes}.  Figure  \ref{fig:RusanovSchemes} depicts a range of potential wave speed estimates in the $x$-$t$ plane, associated characteristic lines and well-known numerical schemes.  The different wave paths $x/t=\hat s$ may be viewed as emerging from the solution of the Riemann  problem for the linear advection equation  (\ref{LAE-1}) with characteristic speed  $\hat s(\beta)$, as suggested in \cite{Toro:2001a},
\cite{Toro:2024a}.

Rusanov-type fluxes for the linear advection equation (\ref{LAE-1})  result from (\ref{LAE-3})  and   (\ref{LAE-4}), namely
\begin{equation}                                    \label{LAE-6}
      f^{Rus}_{i+\frac{1}{2}} = \frac{1}{2}(1 + \beta) f(q^{n}_{i+1}) + \frac{1}{2}(1 - \beta) f(q^{n}_{i}) \;,
\end{equation}
which substituted into the conservative formula (\ref{LAE-2}) yields the three-point scheme
\begin{equation}                                     \label{LAE-7}
    q^{n+1}_{i} = b_{-1}q^{n}_{i-1}+ b_{0}q^{n}_{i}+b_{1}q^{n}_{i+1}\;,
\end{equation}
with coefficients 
\begin{equation}                                     \label{LAE-8}
    b_{-1}=\frac{1}{2}c(1+\beta)\;, \hspace{2mm} b_{0}=(1-\beta c)\;, \hspace{2mm} b_{1}=\frac{1}{2}c(-1+\beta)\;
\end{equation}
and Courant number $c$
\begin{equation}                                      \label{LAE-9}
      c = \frac{\lambda \Delta t} {\Delta x}  \;.
\end{equation}
Assumption (\ref{LAE-4}) gives rise to a family of Rusanov-type schemes (\ref{LAE-7})-(\ref{LAE-8}). As indicated, these may be interpreted as emerging from solving the  Riemann problem for (\ref{LAE-1}) with speed $ \beta \lambda$.  Depending on whether $\beta< 1$ or $\beta > 1$  in (\ref{LAE-4}),  the new characteristic line results from  {\bf a pull ($\beta< 1$)} or {\bf a push ($\beta> 1$)} of the true characteristic line (Godunov's method GodU),  as illustrated in Figure \ref{fig:RusanovSchemes}.   A pull ($\beta< 1$)  results in an underestimate of the true wave speed and  a push ($\beta> 1$)  results in an overestimate of the true wave speed.  In both cases an error is inherited.  Particular values of $\beta(c)$ reproduce some well-known schemes in the literature \cite{Toro:2009a}, as seen in Table \ref{Table:betaschemes}.
\begin{center}
\begin{table}[h]
\centering
\begin{tabular}{|c|c|c|c|c|c|c|} \hline  
Scheme         &$b_{-1}$ & $b_{0}$ & $b_{1}$ & $\beta$ & Monotone & $c_{lim}$ \\ \hline
Lax-Friedrichs &$\frac{1}{2}(1+c)$ & 0 & $\frac{1}{2}(1-c)$ & $\beta_{LF}=1/c$ & Yes & 1                \\ \hline 
Lax-Wendroff   &$\frac{1}{2}(1+c)c$ & $1-c^{2}$ &$-\frac{1}{2}(1-c)c$ & $\beta_{LW}=c$ & No & 1 \\ \hline  
FORCE          &$\frac{1}{4}(1+c)^{2}$ & $\frac{1}{2}(1+c^{2})$ & $\frac{1}{4}(1-c)^{2}$ & $\beta_{FO}=\frac{(1+c^{2})}{2c}$ &  Yes & 1  \\ \hline
Godunov upwind          &$c$ & $1-c$ & $0$ & $\beta_{GU}=1$ &  Yes & 1  \\ \hline
Godunov centred   &$\frac{1}{2}(1+2c)c$ & $1-2c^{2}$ &$-\frac{1}{2}(1-2c)c$ & $\beta_{GC}=2c$ & No & $\frac{1}{2}\sqrt 2$   \\ \hline 
FTCS   &$\frac{1}{2}$ & $1$ &$-\frac{1}{2}$ & $\beta_{CS}=0$ & No & *   \\ \hline 
\end{tabular} 
\vspace{0.2cm}
\caption{Rusanov-type three-point schemes.  Coefficients $b_{k}$, wave speed coefficient $\beta$, monotonicity and linear stability limit $c_{lim}$.  Symbol * denotes unconditional instability.  For an interpretation of the schemes see Figure \ref{fig:RusanovSchemes}.  For background on all the methods listed see \cite{Toro:2009a}.}      
\label{Table:betaschemes}
\end{table}
\end{center}
\begin{figure}[h]
      \centerline{
      \includegraphics[scale=0.8]{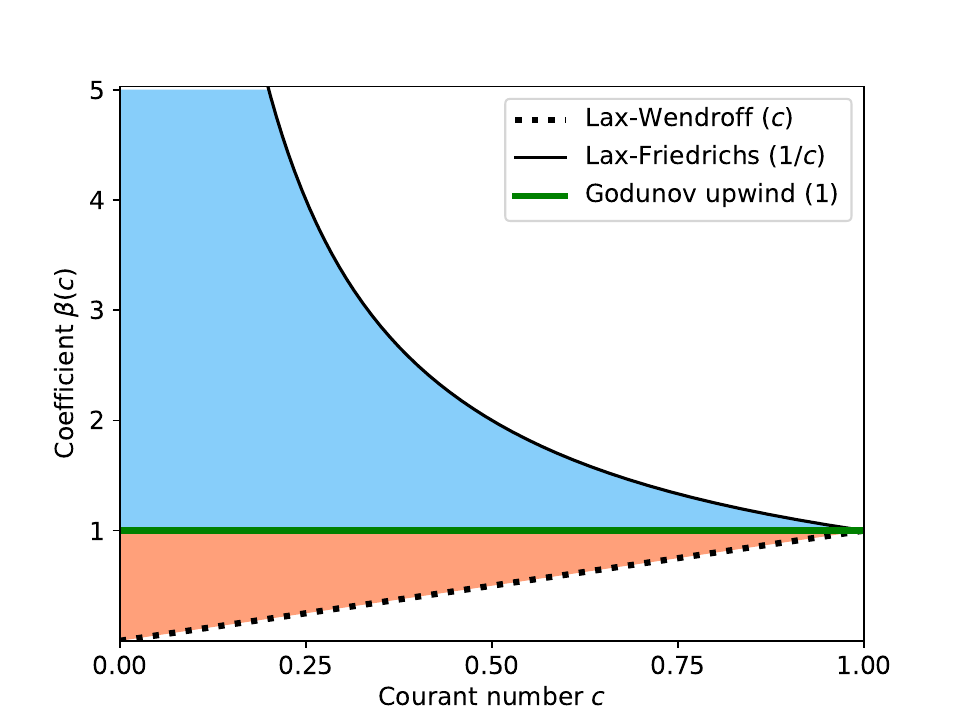}
      }
      \caption{The stability region of the class of Rusanov-type schemes lies between the Lax-Wendroff  ($\beta=c$) and the Lax-Friedrichs schemes ($\beta=1/c$),  and is defined by $c \le \beta \le 1/c$.  The region of monotone schemes (blue) lies between the Godunov upwind ($\beta=1$) and the Lax-Friedrichs ($\beta=1/c$) schemes,  that is  $1 \le \beta \le 1/c$.  The region of non-monotone  but stable schemes (salmon) lies between Lax-Wendroff ($\beta=c$) and Godunov upwind ($\beta=1$),
      that  is $c \le \beta < 1$.  Note that the region of monotone schemes (blue) extends to $\infty$ as $c$ tends to zero.  The rest of the rectangular region (white) identifies unconditionally unstable schemes. }
      \label{fig:FullStability}
\end{figure}
Next we analyse the Rusanov-type schemes emerging from (\ref{LAE-7})-(\ref{LAE-8}), in terms of monotonicity, linear stability and local truncation error.

\subsection{Properties of a class of Rusanov-type schemes}

There are two major properties of schemes that are of interest here,  namely monotonicity and stability.  Local truncation error analysis is also useful. \\

\noindent{\bf Monotonicity. }  Monotonicity of the linear scheme (\ref{LAE-7})-(\ref{LAE-8}) requires non-negativity of coefficients (\ref{LAE-8}).  Simple algebraic manipulations imply the following two cases:
\begin{equation}                                     \label{eq:propRUS1}
\left. \begin{array}{cc}
   0  < \beta < 1   & \hspace{5mm}  \mbox{Non-monotone schemes}\;,  \\
   1 \le \beta \le \displaystyle \frac{1}{c}  & \hspace{5mm}  \mbox{Monotone schemes}\;. 
   \end{array}\right\}
\end{equation}
Note that an underestimate ($0  < \beta < 1$) results immediately in loss of monotonicity.

\noindent{\bf The viscous form. }  For the purpose of this paper, it is convenient to re-write scheme (\ref{LAE-7})-(\ref{LAE-8}) in viscous form; simple manipulations yield
\begin{equation}                                     \label{LAE-11}
    q^{n+1}_{i} = q^{n}_{i} -\frac{1}{2} c (q^{n}_{i+1}-q^{n}_{i-1}) + \frac{1}{2} d (q^{n}_{i+1}-2q^{n}_{i} + q^{n}_{i-1}) \;,
\end{equation}
with coefficient of viscosity
\begin{equation}                                      \label{LAE-12}
     d =\beta c  \;.
\end{equation}

\noindent{\bf Linear stability. }  As is well known, the linear stability condition for scheme (\ref{LAE-11}) is
\begin{equation}                                      \label{LAE-13}
       c^{2} \le d \le 1 \;, 
\end{equation}
The stability condition (\ref{LAE-13}) implies the following two inequalities 
\begin{equation}                                      \label{LAE-14}
        \beta \ge c  \hspace{3mm}  \mbox{   and   }  \hspace{3mm}  c \le 1/\beta \;, 
\end{equation}
to be discussed shortly.\\

\noindent{\bf Coefficient of numerical viscosity. } Local truncation error analysis of schemes (\ref{LAE-7})-(\ref{LAE-8}) gives the  coefficient of numerical viscosity as
\begin{equation}                                      \label{LAE-15}
      \gamma_{visc} = \frac{1}{2}  \Delta x \lambda \left( \frac{d - c^{2}}{c} \right) = 
      \frac{1}{2}  \Delta x \lambda \left( \beta-c \right)\;.
\end{equation}

We now turn to the two stability conditions (\ref{LAE-14}).  They determine  the stability region displayed in Figure \ref{fig:FullStability}, which lies between the Lax-Wendroff  ($\beta(c)=c$) and the Lax-Friedrichs ($\beta(c)=1/c$) schemes,  defined by $c \le \beta(c) \le 1/c$.  There are two subregions of stability defined by (\ref{LAE-14}).
The region of monotone schemes (blue) lies between the Godunov upwind ($\beta(c)=1$) and the Lax-Friedrichs 
($\beta(c)=1/c$) schemes,  that is  $1 \le \beta(c) \le 1/c$.   Note that the region of monotone schemes (blue) extends to $\infty$ as $c$ tends to zero.   The region of non-monotone  but stable schemes (salmon) between Lax-Wendroff ($\beta(c)=c$) and Godunov upwind ($\beta(c)=1$) is defined by $c \le \beta(c) < 1$.   Schemes outside the displayed stability region are obviously unstable.

We note that, classically,  the analysis of schemes is performed in terms of $d$, the coefficient of viscosity in the viscous form of the scheme (\ref{LAE-11}).  In the present study we formulate the schemes in terms of $\beta(c)$. The benefit of this choice is that the analysis relates directly to wave speed estimates, which is the most prominent ingredient of HLL-type schemes (e.g. Rusanov) in practice.  Moreover, wave speed estimates are the underlying subject of interest to various approaches to the design of numerical methods,  see \cite{Guermond:2016a} and \cite{Toro:2020c}.

\subsection{Classical schemes in the  Rusanov framework}
\label{sec:ClassicalSchemes}

\begin{figure}
      \centerline{
      \includegraphics[scale=0.8,angle=0]{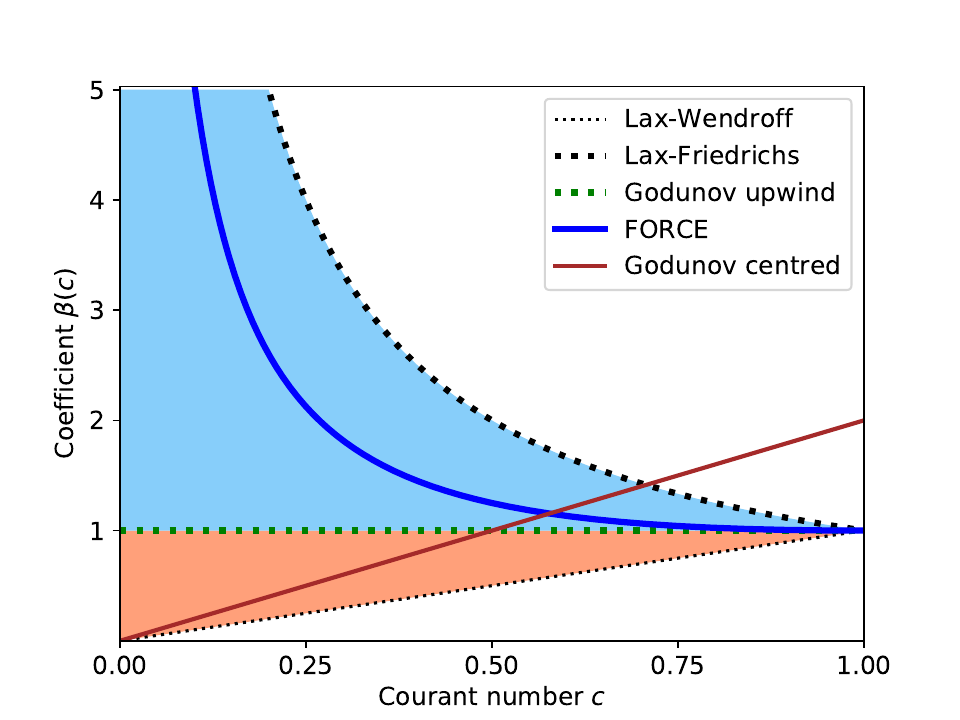}
      }
\caption{Beta functions for conventional methods represented as Rusanov-type schemes. The region of stable schemes is bounded below by the Lax-Wendroff line $\beta(c)=c$  and bounded above by the  Lax-Friedrichs curve $\beta(c)= 1/c$.  As also shown in Figure \ref{fig:FullStability}  the full region of stable schemes  is divided into two sub-regions by the horizontal line $\beta=1$ corresponding  to the Godunov upwind scheme.  
Two more linearly stable schemes are the FORCE scheme given by the curve $\beta(c) = (1+c^{2})/2c$ and the Godunov centred scheme is given by the line $\beta(c)=2c$.  The FORCE method is monotone  in the full range of Courant numbers. The  Godunov centred scheme is non-monotone for $c \le 1/2$  but monotone for $1/2 < c \le  \frac{1}{2} \sqrt{2}=c_{lim}$; the upper limit for $c_{lim}$ is given by the intersection of the Godunov centred scheme with the Lax-Friedrichs scheme.}
\label{fig:ClassicalSchemes}
\end{figure}

Well-known schemes in the literature \cite{Toro:2009a} may be set in the framework of Rusanov-type schemes in terms of appropriate expressions for the wave coefficient $\beta(c)$ in (\ref{LAE-4}). The reader is encouraged to verify that the schemes listed in Table \ref{Table:betaschemes} are identically reproduced by the respective $\beta(c)$ functions.  See Figure \ref{fig:ClassicalSchemes} and  compare it with Figure  \ref{fig:FullStability}, which defines  stability regions.  The schemes in Table \ref{Table:betaschemes} are characterised as follows:
\begin{enumerate}
    \item The last scheme in Table \ref{Table:betaschemes},  called Forward in Time Centred in Space (FTCS),  is obtained from $\beta(c) = \beta_{CS}=0$.  As is well known, this method is unconditionally unstable 
    and $\gamma_{visc}<0$ in (\ref{LAE-15}). 
    \item The Lax-Wendroff method  \cite{Lax:1960a} is obtained from the choice $\beta(c) =  \beta_{LW}=c$.  This 
    scheme is second-order accurate in space and time,  non-monotone and $\gamma_{visc}=0$ in (\ref{LAE-15}).  

    \item The Godunov centred scheme \cite{Godunov:1962a} is obtained from the choice $\beta(c) = \beta_{GC} =2c$. This scheme is non-monotone for $0 < c \le \frac{1}{2}$ and monotone for 
    $\frac{1}{2}\le c \le \frac{1}{2}\sqrt{2}=c_{lim}$. The scheme is linearly stable in the range  $0 \le c \le \frac{1}{2}\sqrt{2}=c_{lim}$.  The stability limit $c_{lim}$ results from the intersection of the function $\beta(c) =2c$ with the Lax-Friedrichs  function  $\beta(c)=1/c$.  See Figure  \ref{fig:ClassicalSchemes}.

    \item The Godunov upwind method \cite{Godunov:1959a} is reproduced identically for $\beta(c) = \beta_{GU}=1$.  This scheme has the smallest truncation error in the family of first-order monotone methods. This is easily seen from (\ref{LAE-8}) and (\ref{LAE-15}), as any decrease of $d=\beta c$ will violate monotonicity.
    
    \item The FORCE method \cite{Toro:1996b}, \cite{Toro:2000a}, \cite{Chen:2004a} is reproduced from the choice $\beta(c) = \beta_{FO}=\frac{1+c^{2}}{2c}$. The FORCE scheme is monotone for the full range of Courant numbers $c \in [0,1]$, as illustrated in Figure  \ref{fig:ClassicalSchemes}.
    
    \item The Lax-Friedrichs method \cite{Lax:1960a} is reproduced from the choice $\beta(c) = \beta_{LF}=\frac{1}{c}$.  This scheme has the largest coefficient of numerical viscosity in the class of linearly stable and monotone schemes and constitutes the upper bound for the class of monotone schemes. See  Figures  \ref{fig:FullStability} and \ref{fig:ClassicalSchemes}. 
    
\end{enumerate}

It can easily be verified that, with the exception of the Godunov centred method,  the coefficients $\beta(c)$ listed in Table \ref{Table:betaschemes} satisfy the following inequalities
\begin{equation}                                   \label{LAE-16}
  \beta_{CS}=0 \le  \beta_{LW}=c \le  \beta_{GU}=1 \le  \beta_{FO}=\frac{1+c^{2}}{2c}  \le  \beta_{LF}=\frac{1}{c}  \;.
\end{equation}
The various coefficients $\beta$, written in increasing order in (\ref{LAE-16}),  give rise to a corresponding sequence of increasing local truncation errors, as can be verified from (\ref{LAE-15}). The increasing sequence of wave speeds $\hat s = \beta \lambda$ determines characteristic lines as depicted in Figure \ref{fig:RusanovSchemes}, which shows a {\bf wave map} in the $x$-$t$ plane of possible three-point schemes of the Rusanov type.   
With reference to Figure \ref{fig:RusanovSchemes} the class of unconditionally unstable schemes lies in the wedge {\it 0OLW}; the class of non-monotone but stable schemes lies in the wedge {\it LWOGU}; the class of monotone and stable schemes lies in the wedge {\it GUOLF}  and the class of unconditionally unstable schemes lies in the wedge {\it LFO$\infty$}.

\subsection{FORCE-$\alpha$ schemes in the Rusanov framework}

Apart from the examples depicted in Figure \ref{fig:ClassicalSchemes} and discussed in Section \ref{sec:ClassicalSchemes},  we show here that there are more monotone methods that can be expressed in the framework of Rusanov-type methods.   One example is furnished by the FORCE-$\alpha$ schemes \cite{Toro:2020b}. These monotone schemes  emanate from the  one-dimensional FORCE scheme \cite{Toro:1996b}, \cite{Toro:2000a}, \cite{Chen:2004a},  formulated as multi-dimensional schemes on general meshes \cite{Toro:2009c}, \cite{Dumbser:2010a}.  The FORCE-$\alpha$ numerical flux for the linear advection equation (\ref{LAE-1}) is
\begin{equation}                                 \label{eq:alfa-0}
      f^{\alpha}_{i+\frac{1}{2}} = \frac{1}{2}(1 + r) f(q^{n}_{i}) + \frac{1}{2}(1 - r) f(q^{n}_{i+1}) \;, \hspace{2mm} r= \frac{1+\alpha^{2}c^{2}}{2 \alpha c} \;. 
\end{equation}
Here $\alpha >0$ is a parameter and $c$ is the Courant number as usual.
\begin{figure}[h]
      \centerline{
      \includegraphics[scale=0.8,angle=0]{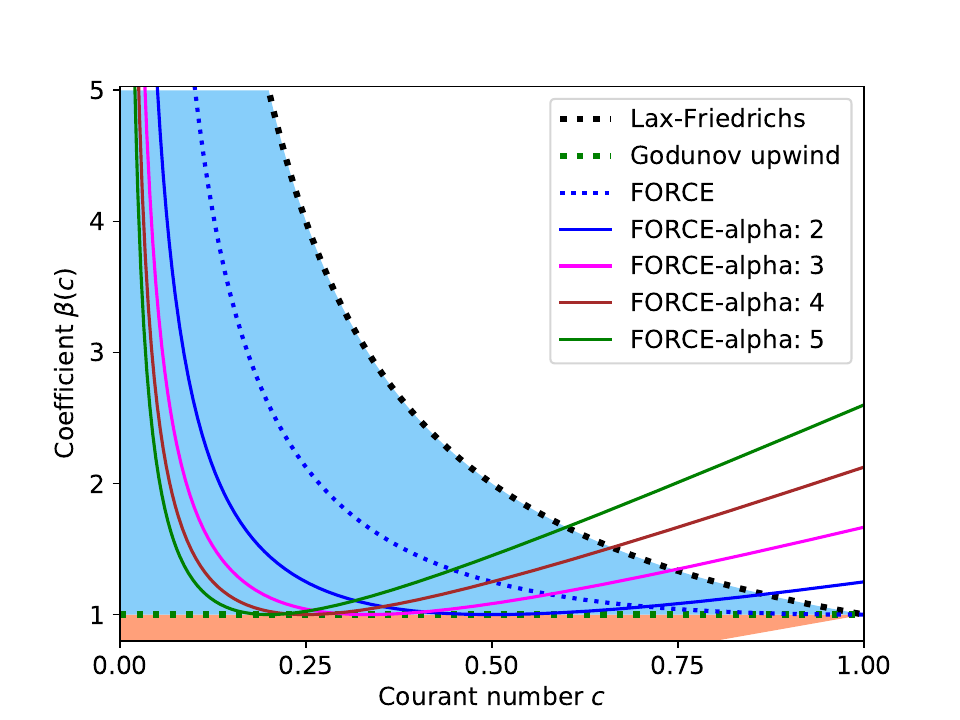}
      }
\caption{$\beta(c;\alpha)$ functions for the FORCE-$\alpha$ schemes, along with three conventional monotone schemes, namely  Godunov upwind,  Lax-Friedrichs and the conventional FORCE method  ($\alpha=1$).   $\beta(c; \alpha)$ functions for the FORCE-$\alpha$ schemes are shown for  $\alpha \in \left\{2,3,4,5 \right\}$.
The stability limit $c_{lim}$ is given by the point of intersection between the $\beta(c;\alpha)$ and the $\beta_{LF}=1/c$ curves.  As $\alpha$ increases $c_{lim}$ decreases. }
\label{fig:AlphaSchemes}
\end{figure}
For the Cartesian-mesh version of the multidimensional FORCE schemes  \cite{Toro:2009c} the specific values 
$\alpha=2$ in two space dimensions  and  $\alpha=3$ in three space dimensions emerge naturally from the multidimensional schemes,  but in the FORCE-$\alpha$ formulation one  regards  
$\alpha$ a free parameter.  From comparing (\ref{eq:alfa-0}) with (\ref{LAE-6}) we obtain the variable coefficient  
$\beta=\beta(c)$ given as
\begin{equation}                        \label{eq:alfa-1}
	\beta= \beta(c; \alpha) = \frac{1}{2}\left(\frac{1}{\alpha c}+\alpha c \right)\;.
\end{equation}
Figure \ref{fig:AlphaSchemes} shows $\beta(c; \alpha)$ functions for the FORCE-$\alpha$ schemes, along with three conventional monotone schemes, namely  Godunov upwind,  Lax-Friedrichs and the conventional FORCE method  ($\alpha=1$);  $\beta(c; \alpha)$  functions for the FORCE-$\alpha$ schemes are shown for 
$\alpha \in \left\{2,3,4,5 \right\}$.  The stability limit $c_{lim}$ of FORCE-$\alpha$ is given by intersection of the corresponding $\beta(c;\alpha)$ curve with that for the Lax-Friedrichs scheme,  that is
\begin{equation}                        \label{eq:alfa-10}
	c_{lim}=\frac{1}{\alpha}\sqrt{2\alpha-1} \;.
\end{equation}
Regarding the stability  limit $c_{lim}$,  the FORCE-$\alpha$ schemes turn out to be quite competitive, for low values of $\alpha$, as compared with other types of Rusanov schemes, as we shall see. This is particularly the case for multidimensional problems,  as discussed in Section  \ref{sec:Rusanov2D}.

\section{Consequences of Perturbing the Exact Speed}
\label{sec:Perturbation}

Here we discuss the consequences of utilising erroneous wave speed estimates by selecting constant 
perturbations of the true wave speed in the Rusanov scheme.  First we discuss monotonicity and stability theoretically and then we illustrate these concepts through some selected numerical experiments.

\subsection{Monotonicity and stability}
\label{sec:MonotonicityStability}

\begin{figure}[t]
      \centerline{
      \includegraphics[scale=0.8,angle=0]{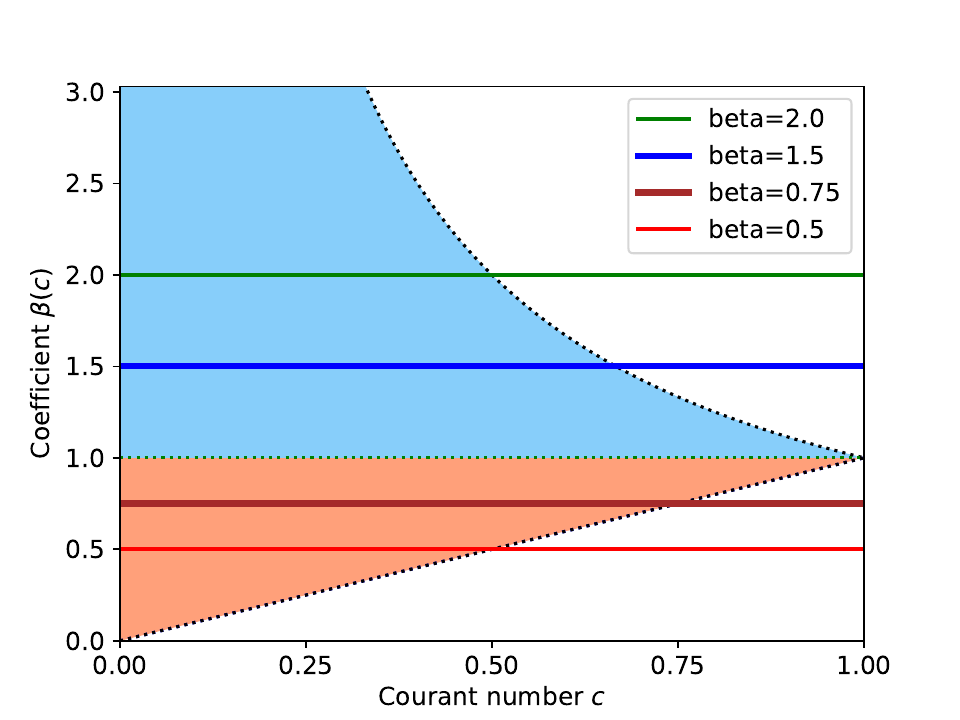}
      }
\caption{Rusanov-type schemes resulting from constant perturbations of the Godunov upwind scheme  ($\beta=1$) of the form
$\beta = 1-\epsilon_{B}$ (underestimate) and $\beta = 1+\epsilon_{T}$ (overestimate).
Two examples shown correspond to underestimation of $\hat s$,  namely $\epsilon_{B}=0.5$ ($\beta=0.5$) and $\epsilon_{B}=0.25$ ($\beta=0.75$).   
Two more examples correspond to overestimation of $\hat s$,   namely  $\epsilon_{T}=0.5$ ($\beta=1.5$) and $\epsilon_{T}=1.0$ ($\beta=2.0$).}
\label{fig:BetaConstant}
\end{figure}

A convenient way of representing Rusanov-type schemes in the monotonicity and stability  regions depicted in Figure  \ref{fig:FullStability} is by considering $\beta(c)$ functions of the form
\begin{equation}                        \label{eq:BetaConst1}
         \beta(c) = \mbox{constant} \in [0, \infty) \;, \hspace{2mm} c \in [0,1] \;.
\end{equation}
See Fig.  \ref{fig:BetaConstant}.
In what follows the discussion regarding $\beta(c)$ is performed in terms of the two parameters $\epsilon_{B}$ 
and $\epsilon_{T}$  introduced in (\ref{LAE-5}).  Beta functions of the form (\ref{eq:BetaConst1}) satisfying (\ref{LAE-5}) offer a realistic approach for assessing the practical scenarios arising in the estimation of the speed $\hat s$.  In reality,  at each interface there are two possible options, namely underestimation (pull,  non-monotone schemes) for $ 0\le \beta(c) <1$ and overestimation (push, monotone schemes) for $ 1 < \beta(c) \le 1/c$.   Figure \ref{fig:BetaConstant}  shows four examples, two in the non-monotone region ($\beta=0.5$, $\beta=0.75)$ and two in the monotone region ($\beta=1.5$, $\beta=2.0)$.

Note that in both cases $\beta=1-\epsilon_{B}$ (non-monotone schemes) and $\beta=1+\epsilon_{T}$ (monotone schemes)  the corresponding stability limit $c_{lim}$ is smaller than the optimal value $c=1$ for the 
Godunov upwind scheme.  That is,  part of the stability range is lost.
In the non-monotone region the upper limit $c_{lim}$ is given by the intersection of  $\beta=1-\epsilon$ with the line $\beta = c$ and hence $c_{lim}=1-\epsilon_{B}$.  As $\epsilon_{B}$ increases $c_{lim}$ decreases. 
In the monotone region the upper limit $c_{lim}$ is given by the intersection of  $\beta=1+\epsilon_{T}$ with the line $\beta = 1/c$ and hence $c_{lim}=\frac{1}{1+\epsilon_{T}}$.  As $\epsilon_{T}$ increases $c_{lim}$ decreases.   As seen in Figure \ref{fig:BetaConstant}  for a given constant displacement $\epsilon$ in $\beta=1-\epsilon$ or $\beta=1+\epsilon$ the stability restriction is more severe in the non-monotone region associated with underestimation of $\hat s$,  than in the monotone region of an overestimated $\hat s$.    This is clearly seen through the following relations
\begin{equation}                                   \label{eq:BetaConst3}
 \mbox{If  }  \epsilon_{B}=\epsilon_{T}=\epsilon  \hspace{3mm} \mbox{  then  }  \hspace{3mm}
 c_{lim}^{under}=1-\epsilon <c_{lim}^{over}= \frac{1}{1+\epsilon} \;.
\end{equation}

\begin{figure}[h]
      \centerline{
      \includegraphics[scale=0.8,angle=0]{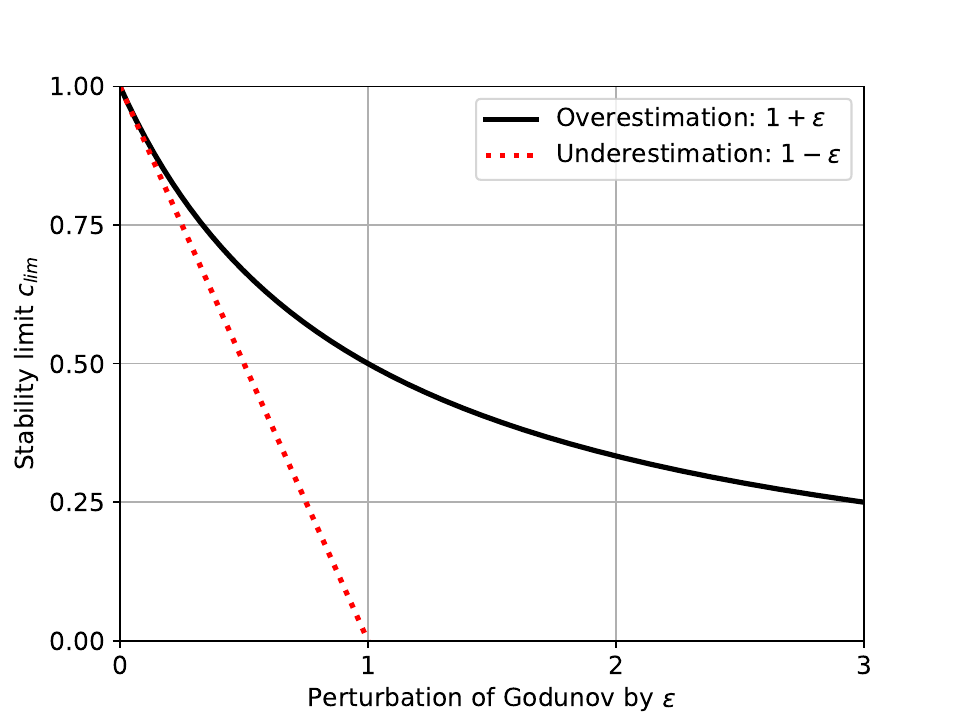}
      }
\caption{Stability limit $c_{lim}$ for Rusanov-type schemes as function of a perturbation $\epsilon$ of the Godunov upwind method.  Underestimation results in non-monotone methods and the stability limit is given by $c_{lim}=1-\epsilon$; when $\epsilon=1$, $c_{lim}=0$.  Overestimation results in monotone methods and the stability limit is given by $c_{lim}=\frac{1}{1+\epsilon}$; when $\epsilon$ tends to $\infty$ $c_{lim}$ tends to zero.}
\label{fig:RusanovStability}
\end{figure}

Figure \ref{fig:RusanovStability} illustrates the stability limits for both under and overestimation. 
Here are some examples.  For $\epsilon_{B}=1/2$ (underestimate),  $c_{lim}^{under}=1/2$, whereas for $\epsilon_{T}=1/2$ (overestimate),  $c_{lim}^{over}=2/3$. This observation supports the aim of seeking upper bounds for the true wave speeds when estimating $\hat s$,  as advocated in  \cite{Guermond:2016a}, \cite{Toro:2020c}.   A conclusion with practical implications is the following.  It is preferable to overestimate the speed  $\hat s$ for two reasons.  Firstly,  the schemes are monotone.  Secondly,   for the same deviation  from the exact speed the schemes have a larger stability range,  as shown by (\ref{eq:BetaConst3}) and illustrated in Figures \ref{fig:BetaConstant} and \ref{fig:RusanovStability}.

\subsection{A numerical example}

\begin{figure}[h]
      \centerline{
      \includegraphics[scale=0.6,angle=0]{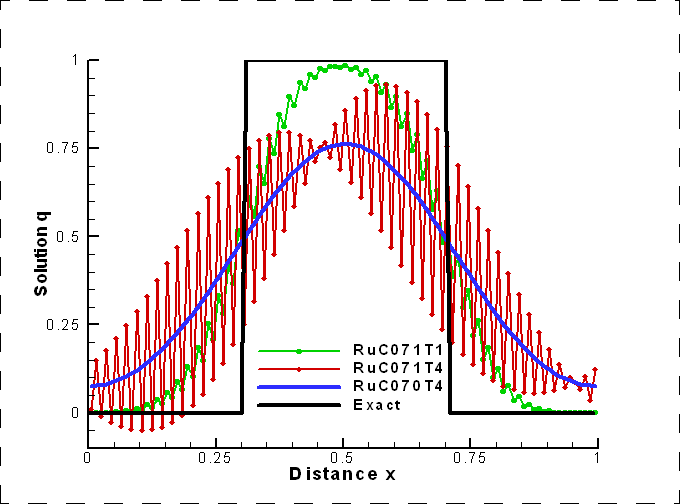}
      }
\caption{Numerical example: $\beta=\sqrt{2}$ gives a stable scheme up to $c_{lim} =\frac{\sqrt{2}}{2} \approx 0.7071$.   Exact and numerical solutions are shown.  
RuC070T4 depicts the result (stable) for $c=0.70<c_{lim}$ at time $T_{out}=4.0$.   
RuC071T1 shows the result (unstable) for $c=0.71>c_{lim}$ at time $T_{out}=1.0$
and RuC071T4 shows the result (more clearly unstable) for $c=0.71>c_{lim}$ at time  $T_{out}=4.0$.}
\label{fig:RusanovExperiments}
\end{figure}

Here we show a numerical example to illustrate the reduction in the stability limit $c_{lim}$ due to deviations  from the true wave speed.  We consider the case of overestimation of $\hat s$ by choosing
$\beta=\sqrt{2}\approx 1.4142136$,  that is $\epsilon_{T}=\sqrt{2}-1 \approx 0.4142136$ in (\ref{LAE-5}). The stability limit for the corresponding Rusanov method is $c_{lim}=1/\sqrt{2}\approx 0.7071068$, well below $c_{lim}=1$ for the Godunov upwind method with $\beta=1$.  To illustrate the situation we solved numerically the following initial-boundary value problem
\begin{equation}                                    \label{eq:BetaConst4}
\left.\begin{array}{cc}
     \mbox{PDE:} & \partial_{t}q + \lambda \partial_{x}q = 0  \;, \hspace{2mm}  x \in [0,1],  \hspace{2mm} t>0 \;,\\
        \mbox{IC:} & q(x,0)= 
        \left\{\begin{array}{ccc}
            0     &   \mbox{if} &    x < 1/4  \;, \\
            1     &   \mbox{if} &     1/4 \le x \le  3/4 \;, \\
            0    &   \mbox{if} &     x > 3/4  \;.
      \end{array}\right.
 \end{array}\right\}
\end{equation}
Here we choose $\lambda=1$ and and impose periodic boundary conditions.  Figure \ref{fig:RusanovExperiments} compares the exact solution  with numerical solutions from Rusanov-type schemes, for two output times: $T_{out}=1$ (T1) and $T_{out}=4$ (T4). Two values for the Courant number are considered: $c=0.70$ just below the stability limit $c_{lim}=0.7071068$ and $c=0.71$ just above the stability limit.  The full blue line shows the numerical solution for $c$ just inside the stable region at time  $T_{out}=4$ (RuC070T4). The solution looks perfectly stable after four cycles within the domain,  in agreement with the theory.  The green line (and dots) shows results for the shorter output time  $T_{out}=1$ (T1),  but for $c$ just outside the stability region (RuC071T1).  In this result the predicted instability is just beginning to manifest.
The red cuve (and dots) (RuC071T4) shows results for $c$ just outside the stability region but for the longer output time  $T_{out}=4$.  The solution is clearly unstable,  as predicted.  The point here is that for any $\beta(c)>1$ above the 
Godunov upwind scheme (overestimate),  the Rusanov method looses part of the stability range,   of the order of  
$30 \%$ in the present example with $\beta=\sqrt{2}$.  As already pointed out, the same phenomenon is observed for $\beta(c)<1$ 
(underestimation), where the solution is in fact non-monotone and the reduction in the stability range is actually more severe than in the overestimation case,   as seen in (\ref{eq:BetaConst3}) and illustrated in Figures \ref{fig:BetaConstant} and \ref{fig:RusanovStability}.

\section{Rusanov-Type Schemes in Two Space Dimensions}
\label{sec:Rusanov2D}

In this section we first formulate Rusanov-type schemes in two space dimensions, based on a simultaneous updating formula,  and then study monotonicity and stability of resulting schemes.

\subsection{Formulation of Rusanov-type schemes} 

We we extend the Rusanov-type schemes as defined previously to the linear advection  equation in two space dimensions on a  Cartesian mesh of size $\Delta x \times \Delta y$
\begin{equation}                        \label{eq:Stab-1}
	\partial_{t} q(x,y,t)+\partial_{x}f(q(x,y,t))+\partial_{y}g(q(x,y,t))=0\;,  \hspace{3mm}
	f(q)=\lambda_{x}q   \;,   \hspace{2mm}  g(q)=\lambda_{y}q \;.
\end{equation}
Without loss of generality we assume $\lambda_{x} \ge 0$ and  $\lambda_{y} \ge 0$.  The corresponding simultaneous update numerical formula is
\begin{equation}                                     \label{eq:Stab-2}
	q_{i,j}^{n+1}=q_{i,j}^n - \frac{\Delta t}{\Delta x}(f^{Rus}_{i+\frac{1}{2},j})-f^{Rus}_{i-\frac{1}{2},j})
    -\frac{\Delta t}{\Delta y}(g^{Rus}_{i, j+\frac{1}{2}} - g^{Rus}_{i, j-\frac{1}{2}}) \;.
\end{equation}
The numerical fluxes in the $x$ and $y$ directions are
\begin{equation}                                     \label{eq:Stab-3}
\left.
\begin{array}{c}
      f^{Rus}_{i+\frac{1}{2},j} = \frac{1}{2}[f(q^{n}_{i,j}) +  f(q^{n}_{i+1,j})]
                                - \frac{1}{2} \hat s_{x} (q^{n}_{i+1,j} - q^{n}_{i,j} ) \;, \\
                                \\
      g^{Rus}_{i,j+\frac{1}{2}} = \frac{1}{2}[g(^{n}_{i,j}) +  g(q^{n}_{i,j+1})]
                                - \frac{1}{2} \hat s_{y} (q^{n}_{i,j+1} - q^{n}_{i,j} ) \;,
      \end{array} \right\}
\end{equation}
where
\begin{equation}                                     \label{eq:Stab-4}
\hat s_x= \beta_{x}  \lambda_{x} \;, \hspace{2mm}
\hat s_y= \beta_{y}  \lambda_{y} \;, \hspace{2mm} 
c_x=\frac{\lambda_x \Delta t}{\Delta x} \;,\hspace{2mm} 
c_y=\frac{\lambda_y \Delta t}{\Delta y} \;.
\end{equation}

\subsection{Monotonicity in two space dimensions}

The Rusanov-type scheme (\ref{eq:Stab-2}) in two space dimensions  with definitions (\ref{eq:Stab-3})-(\ref{eq:Stab-4})  may be written as
\begin{equation}              \label{eq:Stab-5}
    q_{i,j}^{n+1}=\gamma_{-1 }q_{i-1,j}^{n}+\gamma_{0}q_{i,j}^{n} +\gamma_{1}q_{i+1,j}^{n}+
    \delta_{-1}q_{i,j-1}^{n}+\delta_{1}q_{i,j+1}^{n} \;, 
\end{equation}
with coefficients as given in Table \ref{Table:coefficients2D}.  The monotonicity of the linear scheme (\ref{eq:Stab-5}) is determined by the non-negativity of the coefficients in Table \ref{Table:coefficients2D}.  Simple manipulations show that monotonicity,  or otherwise,  is decided exclusively by the coefficients $\gamma_{1}$ and $\delta_{1}$, as the remaining coefficients are all non-negative under the assumptions $\lambda_{x} \ge 0$ and  $\lambda_{y} \ge 0$.   If $\gamma_{1}=\frac{1}{2}\left(-1+\beta_x\right)c_x \ge 0$ then the scheme may be monotone for $\beta_{x}\ge 1$, otherwise it will be non-monotone.  A similar condition applies to 
$\delta_{1}$, that is if $\beta_{y}\ge 1$,  the scheme may be monotone,  otherwise it will be non-monotone.

\tiny
\begin{center}
\begin{table}[h]
\centering
\begin{tabular}{|c|c|c|c|c|} \hline 
 $\gamma_{-1}$            &$\gamma_{0}$             & $\gamma_{1}$  & $\delta_{-1}$ & $\delta_{1}$ \\
\hline
\hline
$\frac{1}{2}\left(1+\beta_x\right)c_x $ & $-\beta_xc_x-\beta_yc_y+1$ & $\frac{1}{2}\left(-1+\beta_x\right)c_x$ & $\frac{1}{2}\left(1+\beta_y\right)c_y$ & $\frac{1}{2}\left(-1+\beta_y\right)c_y$ \\ \hline
\end{tabular} 
\caption{Coefficients $\gamma_{k}$ and $\delta_{l}$ for scheme (\ref{eq:Stab-5}) as functions of Courant numbers 
$c_{x}$, $c_{y}$  and speed coefficients $\beta_{x}$, $\beta_{y}$.}
\label{Table:coefficients2D}
\end{table}
\end{center}
\normalsize

\subsection{Stability in two space dimensions}
	
\begin{figure}
      \centerline{
             \includegraphics[scale=0.4, angle=0]{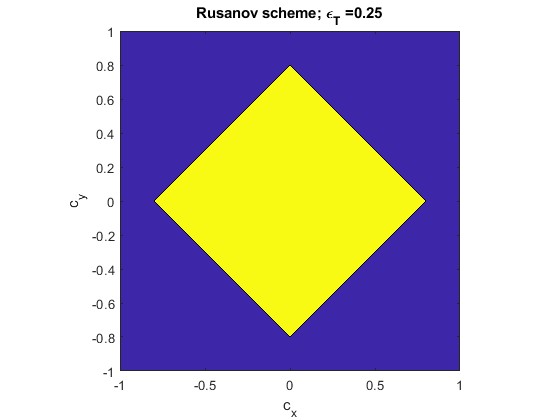}
             \includegraphics[scale=0.4, angle=0]{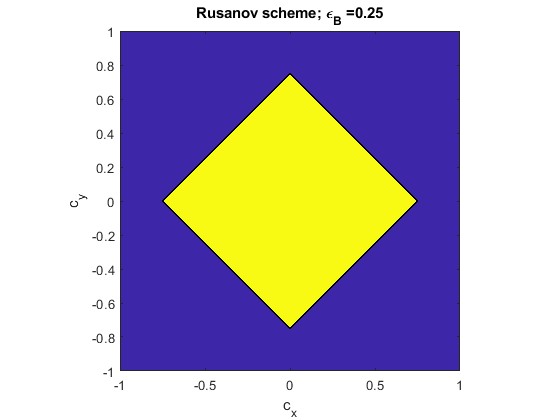}
              }
              \centerline{
              \includegraphics[scale=0.4, angle=0]{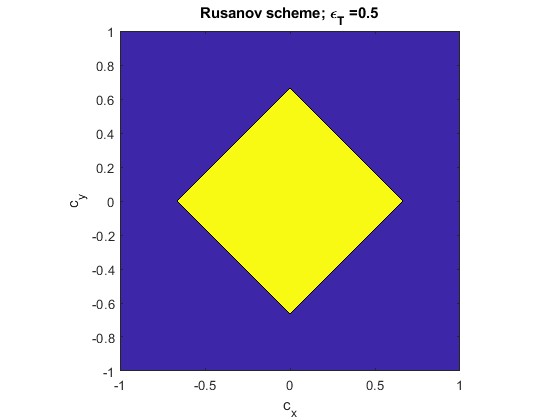}
             \includegraphics[scale=0.4, angle=0]{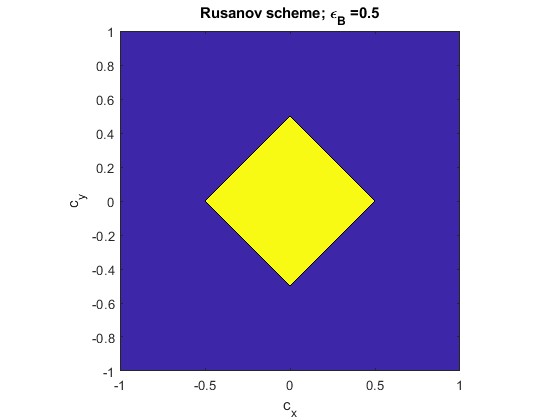}
              }
         \centerline{
             \includegraphics[scale=0.4, angle=0]{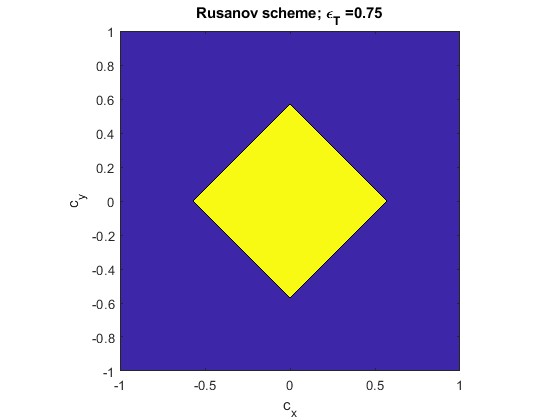}
              \includegraphics[scale=0.4, angle=0]{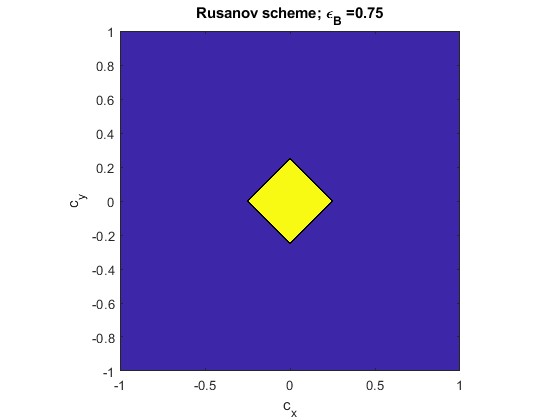}            
      }
      \caption{{\bf Stability regions for two-dimensional  monotone Rusanov schemes}.  The left colum shows 
      stability regions for  the case of overestimated wave speeds,  above $\beta_{GU}=1$.  That is $\beta = 1 +
      \epsilon_{T}$ with  $\epsilon_{T} \in \left\{0.25, 0.5,0.75 \right\}$.  As $\epsilon_{T}$ increases the stability region 
      decreases. The right column shows stability regions for  the case  of under estimated wave speeds,  below $
      \beta_{GU}=1$.  That is $\beta = 1 -\epsilon_{B}$ with $\epsilon_{B} \in \left\{0.25, 0.5,0.75 \right\}$.   As $
      \epsilon_{B}$ increases the stability region decreases.   For equivalent displacements from $\beta_{GU}=1$, the 
      stability region from overestimating (left) is larger than that from underestimating (right) the wave speed $\hat s$.}
      \label{fig:2DstabilityUnderOver}
\end{figure}

\begin{figure}
      \centerline{
              \includegraphics[scale=0.4, angle=0]{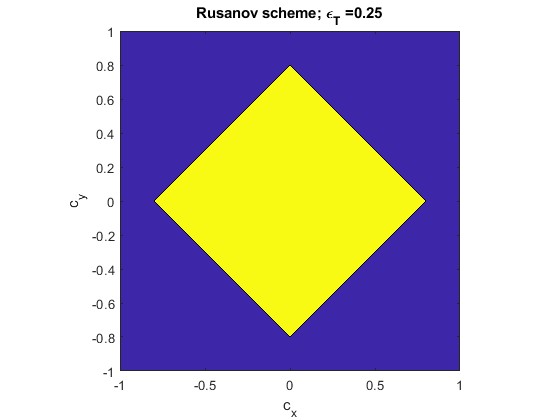}
             \includegraphics[scale=0.4, angle=0]{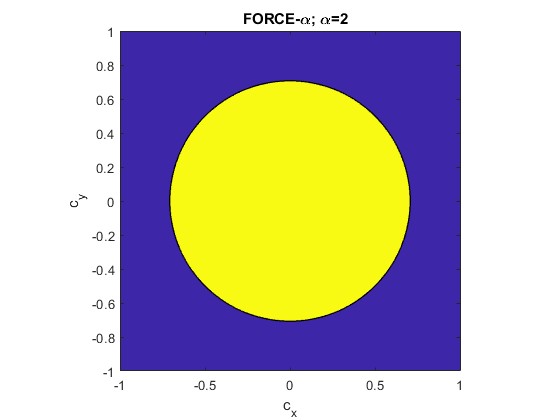}
              }
         \centerline{
                      \includegraphics[scale=0.4, angle=0]{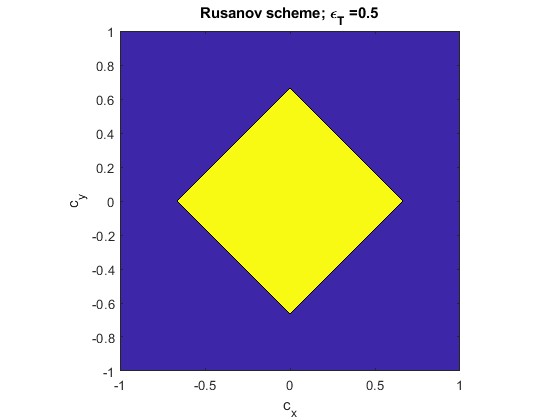}
             \includegraphics[scale=0.4, angle=0]{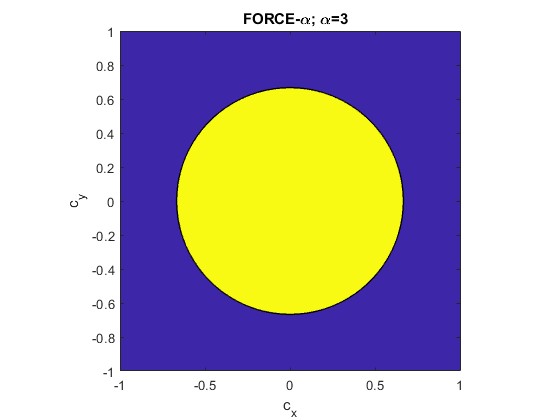}
      }
               \centerline{
                 \includegraphics[scale=0.4, angle=0]{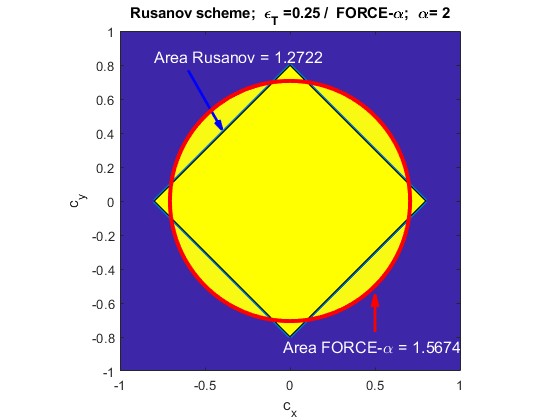}
                \includegraphics[scale=0.4, angle=0]{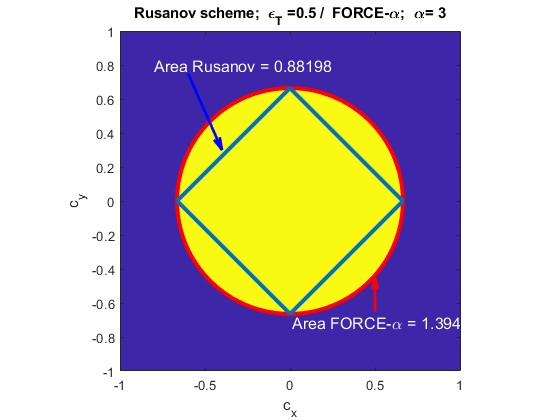}
      }
      \caption{{\bf Stability regions for two-dimensional monotone schemes}.  Top four pictures: comparison of Rusanov 
      monotone   schemes (left column) resulting from overestimation of wave speed,  with $\beta = 1 +\epsilon_{T}$,  
      against the FORCE-$\alpha$ schemes  (right column).  Here $\epsilon_{T} \in \left\{0.25, 0.5 \right\}$ and   $\alpha 
      \in \left\{2, 3 \right\}$.  Bottom two pictures: comparison of size of stability regions for Rusanov-type and 
      FORCE-$\alpha$ schemes.}
      \label{fig:2DstabilityComparison}
\end{figure}

Here we analyse the stability of the  Rusanov-type schemes (\ref{eq:Stab-2}) by means of the von Neumann method \cite{Hirsch:1988a} combined with numerical evaluation of the modulus of the amplification factor for large numbers of parameters,  such as Courant numbers and phase angles.  Results for the stability region for two-dimensional schemes are displayed within a rectangle $[0,  c_{x,lim}]\times[0,c_{y,lim}]$ through two regions,  as depicted in Figures 
\ref{fig:2DstabilityUnderOver} and \ref{fig:2DstabilityComparison},  in which the stable region is shown in yellow. These are obtained by assigning two values: 1  (stability) when for each pair $(c_{x},c_{y})$ for all angles considered the modulus of the amplification factor  is less or equal to $1$.  Otherwise the value $0$ is assigned (instability).  The methodology is fully described in Appendix B of \cite{Toro:2000a} and in Chapter 16 of \cite{Toro:2009a}; see also \cite{Toro:2020b}. 

For the present study we have assumed $\lambda_{x}=\lambda_{y}=1$ and $\beta_{x}=\beta_{y}=\beta: \mbox{constant}$,  for each case.   An entirely analogous procedure applies to three space dimensions.  Stability results from the present study are shown in Figures 
\ref{fig:2DstabilityUnderOver} and \ref{fig:2DstabilityComparison} for the two-dimensional schemes (\ref{eq:Stab-2}).
The stability regions in yellow in the $c_{x}$-$c_{y}$ plane can be described as follows
\begin{equation}                                     \label{eq:Stab-last }
             |c_{x}|+|c_{y}| \le c_{lim} \;,   \hspace{3mm}  c_{lim}= min \left\{c_{x,lim},c_{y,lim}\right\} \;.
\end{equation}
Here  $c_{x,lim}$ is the stability limit for the one-dimensional case in the $x$ direction,  which results from the intersection of $\beta$ with either $1/c$ (Lax-Friedrichs) if $\beta \ge 1$  or with $c$ (Lax-Wendroff) if $\beta < 1$. That is
\begin{equation}                                     \label{eq:Stab-last }
             c_{x,lim} = \left\{
             \begin{array}{ccc} 
             1-\epsilon_{B}                                                  & \mbox{if}   &  \beta < 1 \;\\
             \\
             \displaystyle \frac{1}{1+\epsilon_{T}}      & \mbox{if}   &  \beta \ge 1 \;.
             \end{array}
             \right.
\end{equation}
An analogous definition applies to $c_{y,lim}$ for the $y$ direction.  Figure \ref{fig:RusanovStability} illustrates the stability limits for both under and overestimation in one space dimension. 

Figure  \ref{fig:2DstabilityUnderOver} depicts stability regions for two-dimensional  Rusanov-type schemes for both overestimates (left column) and underestimates (right column) for the wave speed $\hat s$,  for three deviations from the exact value corresponding to the Godunov upwind scheme $\beta_{GU}=1$.  The results shown on the left column correspond to overestimates of the form $\beta = 1 + \epsilon_{T}$,  with   $\epsilon_{T} \in \left\{0.25, 0.5,0.75 \right\}$.    As $\epsilon_{T}$ increases the stability region  decreases (from top to bottom). The results on the right column correspond to underestimated wave speeds of the form $\beta = 1 -\epsilon_{B}$,  with $\epsilon_{B} \in \left\{0.25, 0.5,0.75 \right\}$.   Again,  as  $\epsilon_{B}$ increases the stability region decreases (from top to bottom).   For equivalent displacements from  $\beta_{GU}=1$, the stability region from overestimation (left) is larger than that from underestimation (right).  Recall also that in addition to a reduced stability region,  the results on the right column correspond to non-monotone schemes.  Hence underestimation results in both non-monotone schemes and reduced stability region, as compared with overestimation.  Note also that dependence of the size of the stability region in two space dimensions is entirely consistent with that for the corresponding one-dimensional case; this notion is supported 
by relationship  (\ref{eq:Stab-last }). 

There seems to exist a continuum of monotone schemes,  with no finite limit to $\epsilon_{T}$ in (\ref{LAE-5}) and that,  as $\epsilon_{T}$ tends to $\infty$,  $\beta$ tends to  $\infty$ and the stability limit $c_{lim}$  tends to zero,  and thus 
the area of the  two-dimensional stability region would vanish.  This observation is consistent with the stability regions for the one-dimensional case shown in Figures \ref{fig:BetaConstant} and \ref{fig:RusanovStability}. The stability for the one-dimensional case can be verified along the $c_{x}$ and $c_{y}$ axis in Figure  \ref{fig:2DstabilityUnderOver}; see also Figure \ref{fig:BetaConstant}.  For example,  
for overestimation (left) $c_{x,lim}$ is given by the intersection of the horizontal line $\beta=1+\epsilon_{T}$ with the curve $\beta = \beta_{LF}=1/c$.  
For $\epsilon_{T}=0.25$ $c_{x,lim}=4/5$; 
for $\epsilon_{T}=0.5$ $c_{x,lim}=2/3$;
for $\epsilon_{T}=0.75$ $c_{x,lim}=4/7$.  
Similar observations can be made for the underestimation case.

Figure \ref{fig:2DstabilityComparison} shows stability regions for two classes of  two-dimensional monotone schemes.  The four top pictures show Rusanov-type monotone schemes  for $\beta > 1$ (left column) with $\beta = 1 +\epsilon_{T}$,  compared against the FORCE-$\alpha$ schemes  (right column).  Here  $\epsilon_{T} \in \left\{0.25, 0.5 \right\}$ and $\alpha \in \left\{2, 3 \right\}$.  The bottom two pictures compare the size of the stability regions for both classes of methods.  Note that the FORCE-$\alpha$ schemes are monotone for the full range of Courant numbers within its stability interval,  whereas the Rusanov schemes,  with underestimation of $\beta$,  not only have reduced stability region but are also non-monotone.   

\section{Summary and Concluding Remarks}
\label{sec:Conclusions}

This paper has first addressed the question of providing theoretical wave speed estimates for hyperbolic equations, such as the Euler equations of gas dynamics; a brief review of existing methods and their performance has been given, showing that most methods in use fail to bound the exact wave speeds.  In such scenario, a legitimate question regards the potential consequences of erroneous wave speed estimates on the performance of associated numerical methods.   In the main body of this paper we addressed this question,  first focusing on providing estimates for the single wave speed $\hat s$ that defines the Rusanov flux.   Generally,  one aims at producing overestimates, that is bounds for the true wave speed.  For systems, the subject has not received sufficient attention and most wave-speed estimate schemes in current use,  underestimate the true wave speeds \cite{Toro:2020c}.  On the bases of the linear advection equation in one and two space dimensions we have carried out a detailed analysis of the potential scenarios, paying particular attention to monotonicity and stability properties of the Rusanov-type schemes resulting from wave-speed estimates.  

Two ways of studying wave-speed estimates have been considered.   Firstly,  we have shown that some well-known schemes in the literature may be put in the framework of Rusanov-type schemes,  thus resulting in associated wave speed estimates. These schemes include Lax-Wendroff (non-monotone),  Godunov centred (non-monotone),  FORCE (monotone), FORCE-$\alpha$ (monotone) and Lax-Friedrichs (monotone).   As is well known, all of these methods are linearly stable in one space dimension.  However,  direct extensions of the FORCE and Lax-Friedrichs monotone schemes  via the simultaneous update formulae,  are unstable in two and three space dimensions; this is also true for the Lax-Wendroff method \cite{Toro:2000a}.  Surprisingly,   however, the non-monotone Godunov centred scheme turns out to be linearly stable in two space dimensions.  Secondly,  another class of Rusanov-type schemes was constructed following an approach that resembles what would happen in computational practice.   The schemes result from constant perturbations around the exact speed.

Two situations emerge from constant perturbations.  The first  is {\bf underestimation of the true wave speed},  with $\beta = 1 -\epsilon_{B}$,  for $0 < \epsilon_{B} \le 1$. The resulting schemes are non-monotone,  but are linearly stable in both one and two space dimensions.   As $\epsilon_{B}$ tends to $1$ the stability limit $c_{lim}$ tends to $0$.  These considerations also hold true in two space dimensions; the two-dimensional schemes are non-monotone but stable and the size of the stability region decreases to $0$ as $\epsilon_{B}$ tends to $1$.   The second case  of interest is {\bf overestimation of the true wave speed},  with $\beta = 1+\epsilon_{T}$,  for $0 \le \epsilon_{T} < \infty$.  This case corresponds to monotone and stable schemes in one and two space dimensions.  However, even in one space dimension the stability limit $c_{lim}$ decreases as $\epsilon_{T}$ increases.   
Out of these two situations arising from perturbations of the form $\beta = 1 -\epsilon_{B}$ (underestimation) and 
$\beta = 1+\epsilon_{T}$ (overestimation),  it is  preferable to aim for the latter,  for two reasons: schemes from overestimation are monotone and their stability regions are larger,  as compared with those from underestimation,  for equivalent displacements from the exact speed.  We expect that the conclusions drawn for the two-dimensional case will extend to the three-dimensional case.   
The findings of this paper may prove useful in computational practice. They may also motivate subsequent studies for non-linear hyperbolic systems,  requiring estimates for two or more waves,  such as in HLL and HLLC schemes.

A final remark concerning non-linear equations is in order.  There are indeed theoretical results available regarding monotonicity and entropy stability; see for example \cite{Crandall:1980a},   \cite{Godlewski:1991a}, \cite{Godlewski:2021a}, \cite{Guermond:2024a}.  In such studies it is simply assumed that the maximum wave speeds exist and are available.  In practice one needs to be quantitative and provide estimates for each specific equation or system of equations \cite{Guermond:2016a}, \cite{Toro:2020c}, \cite{Guermond:2024a}.

\begin{center}
       {\bf ACKNOWLEDGEMENTS}
\end{center}
E F Toro acknowledges partial support from project number P130625149,  entitled Development of a Numerical Multiphase Flow Tool for Applications to Petroleum Industry,  funded by Repsol S.A.  Spain.
S A Tokareva was supported by the Laboratory Directed Research and Development Program of Los Alamos National Laboratory.  The LANL unlimited release number is LA-UR-23-32691.

Both authors thank the late Arturo Hidalgo  for hosting a two-week stay in Madrid,  in May of 2023, during which most of the present work was carried out,  with the partial support from the research project PID2020-112517GB-I00 of Ministerio de Econom\'ia y Competitividad (Spain) and by the project CAT200600041M of C\'atedra Fundaci\'on CEPSA.  We are grateful to Professor Lourdes Tello and Professor Arturo Hidalgo  for their hospitality and for organizing social and cultural events that pleasantly supplemented our research activities in Madrid.\\

\noindent{\bf Data Availability.} The datasets generated and analysed during the current study are available from the corresponding author on reasonable request.\\

\noindent{\bf Declarations.} Conflict of interests.  The authors declare that they have no conflict of interests.

\newpage


\newpage

\bibliographystyle{plain}
\bibliography{/Users/toro/Work/BIBrefs/refs-latest}

\end{document}